\documentclass{ws-procs9x6}

\newcommand{\F}{\mathcal{F}}
\newcommand{\R}{\mathcal{R}}

\newcommand{\g}{\mathfrak{g}}

\newcommand{\X}{\mathfrak{X}}

\begin{document}

\title{SOME SUBMANIFOLDS OF ALMOST CONTACT \\
        MANIFOLDS WITH NORDEN METRIC \footnote{Partially supported by Scientific
        researches fund of "St. Cyril and St. Methodius" University of Veliko Tarnovo under contract 
        RD-491-08 from 27. 06. 2008.}}

\author{G.NAKOVA}
\address{Department of Algebra and Geometry, \\ 
Faculty of Pedagogics, University of Veliko Tarnovo, \\
1 Theodosij Tirnovsky Str., 5000 Veliko Tarnovo, Bulgaria\\
E-mail: gnakova@yahoo.com} 
\begin{abstract}
In this paper we study submanifolds of almost contact manifolds with Norden
metric of codimension two with totally real normal spaces. Examples of such
submanifolds as a Lie subgroups are constructed.
\end{abstract}

\keywords{Norden metric, almost contact manifold, submanifold, Lie group}

%%% ----------------------------------------------------------------------
\bodymatter
%%% ----------------------------------------------------------------------

%%%%%%%%%%%%%%%%%%%%%%%%%%%%%%%%%%%%%%%%%%%%%%%%%%%%%%%%%%%%%%%%%%%%%%%%%%%0
\section{Introduction}\label{sec_1}
Let $(M,\varphi,\xi,\eta,g)$ be a $(2n+1)$-dimensional almost contact
manifold with Norden metric, i.~e. $(\varphi,\xi,\eta)$ is an almost 
contact structure \cite{Bler} and $g$ is a metric \cite{GaMGri} on 
$M$ such that
\begin{equation}\label{1}
\begin{array}{l}
\varphi^2X=-id+\eta\otimes\xi, \qquad \eta(\xi)=1, \\ 
g(\varphi X,\varphi Y)=-g(X,Y)+\eta(X)\eta(Y), \\
\end{array}
\end{equation}
where $id$ denotes the identity transformation and $X$, $Y$ are 
differentiable vector fields on $M$, i.~e. $X, Y \in \X(M)$. 
The tensor $\tilde{g}$ given by \\ $\tilde{g} (X,Y)=g(X,\varphi Y)+
\eta(X)\eta(Y)$ is a Norden metric, too. Both metrics $g$ and
$\tilde{g}$ are indefinite of signature $(n+1,n)$.
\par
Let $\nabla$ be the Levi-Civita connection of the metric $g$. The tensor
field $F$ of type $(0,3)$ on $M$ is defined by
\[
F(X,Y,Z)=g((\nabla_X\varphi)Y,Z).
\]
\par
A classification of the almost contact manifolds with Norden metric
with respect to the tensor $F$ is given in \cite{GaMGri} and 
eleven basic classes  \\
$\F_i (i=1,2,\dots ,11)$ are obtained.
\par
Let $R$ be the curvature tensor field of $\nabla$ defined by
\[
    R(X,Y)Z=\nabla_X \nabla_Y Z - \nabla_Y \nabla_X Z -
    \nabla_{[X,Y]}Z.
\]
The corresponding tensor field of type $(0,4)$ is
determined as follows
\[
    R(X,Y,Z,W)=g(R(X,Y)Z,W).
\]
\par
Let $(\overline {M}, \overline {\varphi }, \overline {\xi }, \overline {\eta }, g)$
$(dim \overline {M}=2n+3)$ be an almost contact manifold with Norden metric and let
$M$ be a submanifold of $\overline {M}$. Then for each point $p \in
\overline {M}$ we have 
\[
T_p\overline {M}=T_pM\oplus (T_pM)^\bot , 
\]
where $T_pM$ and $(T_pM)^\bot $ are the tangent space and the normal space of 
$\overline {M}$ at $p$ respectively. When the submanifold $M$ of $\overline {M}$
is of codimension 2 we denote $(T_pM)^\bot $ by $\alpha =\{N_1, N_2\}$, i.e.
$\alpha $ is a normal section of $M$.
\par
Let $\alpha $ be a 2-dimensional section in $T_p\overline {M}$. Let us recall
a section $\alpha $ is said to be 
\begin{itemize}
    \item {\it non-degenerate}, {\it weakly isotropic}
or {\it strongly isotropic} if the rank of the restriction of the metric $g$ on
$\alpha $ is 2, 1 or 0 respectively;
		\item of {\it pure} or {\it hybrid} type if the restriction of $g$ on 
		$\alpha $ has a signature $(2,0), \, (0,2)$ or $(1,1)$ respectively;
		\item {\it holomorphic} if $\overline {\varphi }\alpha =\alpha $;
		\item $\overline \xi $- section if $\overline \xi \in \alpha $;
		\item {\it totally real} if $\overline {\varphi }\alpha \bot \alpha $.
\end{itemize}
Submanifolds $M$ of $\overline {M}$ of codimension 2 with a non-degenerate
of hybrid type normal section $\alpha $ are studied. In \cite{NakGriI} two
basic types of such submanifolds are considered: $\alpha $ is a holomorphic section
and $\alpha $ is a $\overline \xi $- section. In \cite{NakGriII} the normal section 
$\alpha =\{N_1, N_2\}$ is such that $\overline {\varphi }N_1\notin \alpha $,
$\overline {\varphi }N_2\in \alpha $. In this paper we consider submanifolds $M$
of $\overline {M}$ of codimension 2 in the case when the normal section $\alpha $
is a non-degenerate of hybrid type and $\alpha $ is a totally real. The totally
real sections $\alpha $ are two types: $\alpha $ is non-orthogonal to 
$\overline \xi $ and $\alpha $ is orthogonal to $\overline \xi $.
\section{Submanifolds of codimension 2 of almost contact manifolds
with Norden metric with totally real non-orthogonal to 
$\overline \xi $ normal spaces}\label{sec_2}
Let $(\overline {M}, \overline {\varphi }, \overline {\xi }, \overline {\eta }, g)$
$(dim \overline {M}=2n+3)$ be an almost contact manifold with Norden metric and let
$M$ be a submanifold of codimension 2 of $\overline {M}$. We assume that there
exists a normal section $\alpha =\{N_1, N_2\}$ defined globally over the 
submanifold $M$ such that
\begin{itemize}
    \item $\alpha $ is a non-degenerate of hybrid type, i.e.
\begin{equation}\label{2}    
g(N_1,N_1)=-g(N_2,N_2)=1, \quad g(N_1,N_2)=0;
\end{equation}    
		\item $\alpha $ is a totally real, i.e.
\begin{equation}\label{3}    
g(N_1,\overline {\varphi }N_1)=g(N_2,\overline {\varphi }N_2)=
g(N_1,\overline {\varphi }N_2)=g(N_2,\overline {\varphi }N_1)=0;
\end{equation}
		\item $\alpha $ is a non-orthogonal to $\overline \xi $ $(\overline \xi \notin
		T_pM)$ and $\overline \xi \notin \alpha $.
\end{itemize}    
Then we obtain the following decomposition for $\overline \xi $, 
$\overline {\varphi }X$, $\overline {\varphi }N_1$, $\overline {\varphi }N_2$ with
respect to $\{N_1, N_2\}$ and $T_pM$
\begin{equation}\label{4}
\begin {array}{llll}
\overline \xi =\xi _0+aN_1+bN_2; \\
\overline {\varphi }X=\varphi X+\eta ^1(X)N_1+\eta ^2(X)N_2, \, \, X\in\chi (M); \\
\overline {\varphi }N_1=\xi _1; \\
\overline {\varphi }N_2=-\xi _2;
\end{array}
\end{equation}
where $\varphi $ denotes a tensor field of type $(1,1)$ on $M$; 
$\xi _0, \xi _1, \xi _2\in\chi (M)$; $\eta ^1$ and $\eta ^2$ are 1-forms
on $M$; $a, b$ are functions on $M$ such that $(a,b)\neq(0,0)$. We denote
the restriction of $g$ on $M$ by the same letter.
\par
Let $a\neq 0$, $\vert a\vert >b$ and $a^2-b^2=k^2$. Taking into account the 
equalities $(1)\div(4)$ we compute
\begin{equation}\label{5}
\eta ^i(X)=g(X,\xi _i), \, (i=0,1,2);
\end{equation}
\begin{equation}\label{6}
g(\varphi X,\varphi Y)=-g(X,Y)+\eta ^0(X)\eta ^0(Y)-\eta ^1(X)\eta ^1(Y)+
\eta ^2(X)\eta ^2(Y);
\end{equation}
\begin{equation}\label{7}
\begin {array}{lll}
\varphi ^2X=-X+\eta ^0(X)\xi _0-\eta ^1(X)\xi _1+\eta ^2(X)\xi _2; \\ 
\eta ^0(\varphi X)=-a\eta ^1(X)+b\eta ^2(X);\\ \eta ^1(\varphi X)=a\eta ^0(X);
\quad \eta ^2(\varphi X)=b\eta ^0(X);
\end{array}
\end{equation}
\begin{equation}\label{8}
\varphi \xi _0=-a\xi _1+b\xi _2;\quad  \varphi \xi _1=a\xi _0;\quad
\varphi \xi _2=b\xi _0;
\end{equation}
\begin{equation}\label{9}
\begin {array}{lll}
g(\xi _0,\xi _0)=1-a^2+b^2; \quad g(\xi _1,\xi _1)=a^2-1; \\ 
g(\xi _2,\xi _2)=1+b^2;\quad 
g(\xi _0,\xi _1)=g(\xi _0,\xi _2)=0; \\ g(\xi _1,\xi _2)=ab;
\end{array}
\end{equation}
for arbitrary $X, Y\in \chi (M)$.
\par
Now we define a vector field $\xi $, an 1-form $\eta $ and a tensor field
$\phi$ of type $(1,1)$ on $M$ by
\begin{equation}\label{10}
\begin {array}{lll}
\xi =\displaystyle{-\frac{b}{k}}.\xi _1+\displaystyle{\frac{a}{k}}.\xi _2; \\ \\
\eta(X)=\displaystyle{-\frac{b}{k}}.\eta^1(X)+
\displaystyle{\frac{a}{k}}.\eta ^2(X), \quad X\in \chi (M); \\ \\
\phi X=\lambda .\varphi ^3X+\mu .\varphi X, \quad X\in \chi (M);
\end{array}
\end{equation}
where
\[
\lambda _1=\displaystyle{\frac{\epsilon}{k(k+1)}}, \quad
\mu _1=\displaystyle{\frac{\epsilon (1+k^2+k)}{k(k+1)}};
\]
\[
\lambda _2=\displaystyle{\frac{\epsilon}{k(k-1)}}, \quad
\mu _2=\displaystyle{\frac{\epsilon (1+k^2-k)}{k(k-1)}};\quad
\epsilon =\pm1.
\]
Further we consider the following cases for $k$:
\par
1) $k^2\neq1\Longleftrightarrow k\neq \pm1$.
In this case $\phi , \xi , \eta $ are given by $(10)$ and $\lambda =\lambda _1,
\, \mu =\mu _1$ or  $\lambda =\lambda _2, \, \mu =\mu _2$.
\par
2) $k=-1$. We obtain $\phi , \xi , \eta $ from $(10)$ by $k=-1$ and 
$\lambda =\lambda _2=\displaystyle{\frac{\epsilon }{2}}, \\ \mu =\mu _2=
\displaystyle{\frac{3\epsilon }{2}}$.
\par
3) $k=1$. We obtain $\phi , \xi , \eta $ from $(10)$ by $k=1$ and 
$\lambda =\lambda _1=\displaystyle{\frac{\epsilon }{2}}, \\ \mu =\mu _1=
\displaystyle{\frac{3\epsilon }{2}}$.
\par
Using $(5)\div(10)$ we verify that $(\phi , \xi , \eta )$ is an almost
contact structure on $M$ and the restriction of $g$ on $M$ is Norden metric.
Thus, the submanifolds $(M, \phi , \xi , \eta , g)$ of $\overline M$ 
considered in 1), 2), 3)
are $(2n+1)$-dimensional almost contact manifolds with Norden metric.
\par
Denoting by $\overline \nabla $ and $\nabla $ the Levi-Civita connections
of the metric $g$ in $\overline M$ and $M$ respectively, the formulas of
Gauss and Weingarten are
\begin{equation}
\begin {array}{lll}\label{11}
\overline \nabla _XY=\nabla _XY+g(A_{N_1}X,Y)N_1-g(A_{N_2}X,Y)N_2; \\ \\
\overline \nabla _XN_1=-A_{N_1}X+\gamma (X)N_2; \\ \\
\overline \nabla _XN_2=-A_{N_2}X+\gamma (X)N_1, \quad X,Y \in \chi (M);\\
\end{array}
\end{equation}
where $A_{N_i} \, (i=1,2)$ are the second fundamental tensors and $\gamma $
is an 1-form on $M$.
\section{Submanifolds of codimension 2 of almost contact manifolds
with Norden metric with totally real orthogonal to 
$\overline \xi $ normal spaces}\label{sec_3}
Let $(\overline {M}, \overline {\varphi }, \overline {\xi }, \overline {\eta }, g)$
$(dim \overline {M}=2n+3)$ be an almost contact manifold with Norden metric and let
$M$ be a submanifold of codimension 2 of $\overline {M}$. We assume that there
exists a normal section $\alpha =\{N_1, N_2\}$ defined globally over the 
submanifold $M$ such that
\begin{itemize}
    \item $\alpha $ is a non-degenerate of hybrid type, i.e.
    the equality $(2)$ holds;
		\item $\alpha $ is a totally real; 
		\item $\alpha $ is orthogonal to $\overline \xi $, i.e. $\overline \xi \in
		T_pM$.
\end{itemize}    
Then from $(4)$ by $a=b=0$ we obtain the following decomposition with respect to 
$\{N_1, N_2\}$ and $T_pM$
\begin{equation}\label{12}
\begin {array}{llll}
\overline \xi =\xi _0; \\
\overline {\varphi }X=\varphi X+\eta ^1(X)N_1+\eta ^2(X)N_2, \, \, X\in\chi (M); \\
\overline {\varphi }N_1=\xi _1; \\
\overline {\varphi }N_2=-\xi _2.
\end{array}
\end{equation}
Substituting $a=b=0$ in $(7), (8), (9)$ we have
\begin{equation}\label{13}
\begin {array}{llll}
\eta ^0(\varphi X)=\eta ^1(\varphi X)=\eta ^2(\varphi X)=0; \\
\varphi \xi _0=\varphi \xi _1=\varphi \xi _2=0; \\
g(\xi _0,\xi _0)=g(\xi _2,\xi _2)=1; \, g(\xi _1,\xi _1)=-1; \\
g(\xi _0,\xi _1)=g(\xi _0,\xi _2)=g(\xi _1,\xi _2)=0.
\end{array}
\end{equation}
Now we define a vector field $\xi $, an 1-form $\eta $ and a tensor field
$\phi$ of type $(1,1)$ on $M$ by 
\begin{equation}\label{14}
\begin {array}{llll}
\xi =t_0\xi _0-t_2\xi _2; \\
\eta (X)=t_0\eta ^0(X)-t_2\eta ^2(X), \quad X\in\chi(M); \\
\phi X=\varphi X+t_0\{\eta ^1(X).\xi _2+\eta ^2(X).\xi _1\}+ \\ 
t_2\{\eta ^0(X).\xi _1+\eta ^1(X).\xi _0\};
\end{array}
\end{equation}
where $t_0, \, t_2$ are functions on $M$ and $t_0^2+t_2^2=1$.
\par
Using $(5), (6), (13), (14)$ we verify that $(\phi , \xi , \eta )$ is an almost
contact structure on $M$ and the restriction of $g$ on $M$ is Norden metric.
So, the submanifolds $(M, \phi , \xi , \eta , g)$ of $\overline M$
are $(2n+1)$-dimensional almost contact manifolds with Norden metric.
The formulas of Gauss and Weingarten are the same as those in section 2.
\section{Examples of submanifolds of codimension 2 of almost contact manifolds
with Norden metric with totally real normal spaces}\label{sec_4}
In \cite{GN-wshop2006} a Lie group as a 5-dimensional almost contact
manifold with Norden metric of the class $\F_9$ is constructed. We will use
this Lie group to obtain examples of submanifolds considered in sections 2 and 3.
\par
First we recall some facts from \cite{GN-wshop2006} which we need. Let 
$\g$ be a real Lie algebra with a global basis of left invariant vector
fields $\{X_1,X_2,X_3,X_4,X_5\}$ and $G$ be the associated with $\g$ real 
connected Lie group. The almost contact structure 
$(\overline \varphi , \overline \xi , \overline \eta )$
and the Norden metric $g$ on ${G}$ are defined by:
\begin{equation}\label{15}
\begin{array}{llll}
\varphi X_i=X_{2+i}, \, \, \, \varphi X_{2+i}=-X_i, \, \, \, \varphi X_5=0, \, \, \, 
(i=1,2); \\
g(X_i,X_i)=-g(X_{2+i},X_{2+i})=g(X_5,X_5)=1, \, \, \, (i=1,2); \\ [4pt]
g(X_j,X_k)=0, \quad  (j\neq k, \, \, \, j, k=1,2,3,4,5); \\
\overline \xi =X_5, \quad \overline {\eta }(X_i)=g(X_i,X_5), \, \, (i=1,2,3,4,5).
\end{array}
\end{equation}
The commutators of the basis vector fields are given by:
\begin{equation}\label{16}
\begin{array}{lll}
[X_1,X_2]=-[X_1,X_3]=aX_4, \quad [X_2,X_3]=aX_2+aX_3, \cr \cr
[X_3,X_4]=-[X_2,X_4]=aX_1, \quad [X_2,\overline \xi ]=2mX_1, \cr \cr
[X_3,\overline \xi ]=-2mX_4, \quad [X_1,X_4]=[X_1,\overline \xi ]=[X_4,\overline \xi ]=0, 
\cr a,m \in {\R}.
\end{array}
\end{equation}
So, the manifold $(G, \overline \varphi , \overline \xi , \overline \eta , g)$ is an
almost contact manifold with Norden metric in the class $\F_9$.
\begin{theorem}\label{4.1}
\cite{Warner} Let ${G}$ be a Lie group with a Lie algebra ${\g}$ and 
$\widetilde b$
be a subalgebra of ${\g}$. There exists an unique connected Lie subgroup $H$ of
${G}$ such that the Lie algebra $b$ of $H$ coincides with $\widetilde b$.
\end{theorem}
\par
From the equalities for the commutators of the basis vector fields 
$\{X_1,X_2,X_3,X_4,X_5\}$ it follows that the 3-dimensional subspaces of ${\g}$
$\bf b_1$ with a basis $\{X_1,X_2,X_3\}$,
$\bf b_2$ with a basis $\{X_1,X_3,X_4\}$ and
$\bf b_3$ with a basis $\{X_1,X_4,\overline \xi \, \}$
are closed under the bracket operation. Hence $b_i$ \, $(i=1,2,3)$ are
real subalgebras of ${\g}$. Taking into account Theorem 4.1 
we have there exist Lie subgroups $H_i \, \,
(i=1,2,3)$ of the Lie group $G$ with Lie algebras $b_i$ \, $(i=1,2,3)$
respectively. 
The normal spaces $\alpha _i \, \, (i=1,2,3)$ of the submanifolds
$H_i \, \, (i=1,2,3)$ of $G$ are:
$\alpha _1=\{X_4,\overline \xi \}$,  
$\alpha _2=\{X_2,\overline \xi \}$, 
$\alpha _3=\{X_2,X_3\}$. Because of $(15)$ we have $\alpha _1$ is $\overline \xi $-
section of hybrid type, $\alpha _2$ is $\overline \xi $- section of pure type and
$\alpha _3$ is a totally real orthogonal to $\overline \xi $ section of hybrid type.
So, the submanifold $H_3$ of $G$ is of the same type submanifolds considered in section 3.   
\par
We choose the unit normal fields of $H_3$ $N_1=X_2$ and $N_2=X_3$. For an arbitrary
$X\in\chi(H_3)$ we have $X=x^1X_1+x^4X_4+\overline {\eta }(X)\overline \xi$. Taking 
into account $(15)$ we compute
\begin{equation}\label{17}
\begin{array}{llll}
\overline \xi=\xi _0; \\
\overline {\varphi }X=-x^4X_2+x^1X_3; \\
\overline {\varphi }X_2=X_4; \\
\overline {\varphi }X_3=-X_1.
\end{array}
\end{equation}
From $(12), (17)$ it follows
\begin{equation}\label{18}
\begin{array}{lll}
\eta ^0(X)=\overline {\eta }(X); \\
\varphi X=0; \, \, \eta ^1(X)=-x^4; \, \, \eta ^2(X)=x^1; \\
\xi _1=X_4; \, \, \xi _2=X_1.
\end{array}
\end{equation}
Substituting $(18)$ in $(14)$ for the almost contact structure on $H_3$ we
obtain
\begin{equation}\label{19}
\begin {array}{llll}
\xi =t_0\overline \xi -t_2X_1; \\
\eta (X)=t_0\overline {\eta }(X)-t_2x^1; \\
\phi X=t_0\{-x^4X_1+x^1X_4\}+ 
t_2\{\overline {\eta }(X)X_4-x^4\overline \xi \};
\end{array}
\end{equation}
where $t_0, \, t_2 \in {\R}$ and $t_0^2+t_2^2=1$.
\par
Using the well known condition for the Levi-Civita connection $\nabla $ of
$g$
\begin{equation}\label{20} 
\begin{array}{l}
2g(\nabla _XY,Z)=Xg(Y,Z)+Yg(X,Z)-Zg(X,Y)+g([X,Y],Z) \\ [4pt]
\phantom{2g(\nabla _XY,Z)=}
+g([Z,X],Y)+g([Z,Y],X) \\
\end{array}
\end{equation}
we get the following equation for the tensor $F$ of $H_3$
\begin{equation}\label{21} 
\begin{array}{l}
F(X,Y,Z)=\displaystyle{\frac{1}{2}}\left\{g([X,\phi Y]-\phi [X,Y],Z)
\right. \\[4pt]
\phantom{F(X,Y,Z)=}
\left.
+g(\phi [Z,X]-[\phi Z,X],Y) \right. \\ [4pt]
\left.
+g([Z,\phi Y]-[\phi Z,Y],X)\right\}, \quad X, Y, Z \in\chi(H_3). 
\end{array}
\end{equation}
From $(16)$ we have $[X_1,X_4]=[X_1,\overline \xi ]=[X_4,\overline \xi ]=0$.
Having in mind the last equalities, $(19)$ and $(21)$ for the tensor $F$ of
$H_3$ we obtain $F=0$. Thus, the submanifold $(H_3, \phi , \xi , \eta , g)$ 
of $G$, where $(\phi , \xi , \eta )$ is defined by (19)
is an almost contact manifold with Norden metric in the class ${\F}_0$.
\par
In order to construct an example for a submanifold from section 2 we 
make the following change of the basis of ${\g}$
\begin{equation}\label{22}
\left(
\begin{array}{l}
E_1 \cr E_2\cr E_3 \cr E_4 \cr E_5
\end{array}
\right)=T^T
\left(
\begin{array}{l}
X_1\cr X_2\cr \overline \xi \cr X_3 \cr X_4
\end{array}
\right), \quad T=
\left(
\begin{array}{lccll}
1 & 0 & 0 & 0 & 0 \cr 
0 & \frac{\sqrt{3}}{2} & \frac{1}{2} & 0 & 0 \cr 
0 & -\frac{1}{2} & \frac{\sqrt{3}}{2} & 0 & 0 \cr
0 & 0 & 0 & 1 & 0 \cr 
0 & 0 & 0 & 0 & 1
\end{array}
\right)\in O(3,2).
\end{equation} 
Taking into account $(16)$ and $(22)$ we compute the commutators of the basis
vector fields $\{E_1,E_2,E_3,E_4,E_5\}$ of ${\g}$
\begin{equation}\label{23}
\begin{array}{lll}
[E_1,E_2]=\displaystyle{\frac{\sqrt{3}}{2}}aE_5, \quad 
[E_1,E_3]=\displaystyle{\frac{1}{2}}aE_5,\quad [E_3,E_5]=
-\displaystyle{\frac{1}{2}}aE_1, \cr \cr
[E_2,E_5]=-\displaystyle{\frac{\sqrt{3}}{2}}aE_1, \quad
[E_1,E_4]=-aE_5, \cr \cr
[E_2,E_4]=\displaystyle{\frac{3}{4}}aE_2+\displaystyle{\frac{\sqrt{3}}{4}}aE_3+
\displaystyle{\frac{\sqrt{3}}{2}}aE_4-mE_5, \cr \cr
[E_3,E_4]=\displaystyle{\frac{\sqrt{3}}{4}}aE_2+\displaystyle{\frac{1}{4}}aE_3+
\displaystyle{\frac{1}{2}}aE_4+\sqrt{3}mE_5, \cr \cr
[E_4,E_5]=aE_1, \quad [E_2,E_3]=2mE_1, \quad [E_1,E_5]=0.
\end{array}
\end{equation}
Because of the elements of the matrix $T$ are constants the Jacobi identity for
the vector fields $\{E_1,E_2,E_3,E_4,E_5\}$ is valid. Now, we compute the matrix
$B$ of $\overline \varphi $ and the coordinates of $\overline \xi $ with respect to
the basis $\{E_1,E_2,E_3,E_4,E_5\}$
\begin{equation}\label{24}
B=\left(
\begin{array}{lclrr}
0 & 0 & 0 & -1 & 0 \cr 
0 & 0 & 0 & 0 &  -\frac{\sqrt{3}}{2} \cr 
0 & 0 & 0 & 0 & -\frac{1}{2} \cr
1 & 0 & 0 & 0 & 0 \cr 
0 & \frac{\sqrt{3}}{2} & \frac{1}{2} & 0 & 0
\end{array}
\right) ; \quad \overline \xi =\left(0,-\displaystyle{\frac{1}{2}},
\displaystyle{\frac{\sqrt{3}}{2}},0,0\right).
\end{equation}
From $T\in O(3,2)$ it follows the matrix of the metric $g$ with
respect to the basis $\{E_1,E_2,E_3,E_4,E_5\}$ is the same as the
matrix 
\begin{equation}\label{25}
C=\left(
\begin{array}{lllrr}
1 & 0 & 0 & 0 & 0 \cr 
0 & 1 & 0 & 0 & 0 \cr 
0 & 0 & 1 & 0 & 0 \cr
0 & 0 & 0 & -1 & 0 \cr 
0 & 0 & 0 & 0 & -1
\end{array} 
\right)
\end{equation}
of $g$ with respect to the basis $\{X_1,X_2,\overline \xi ,X_3,X_4\}$.
\par
Using $(23)$ we have that the 3-dimensional subspace $b$ of ${\g}$ with
a basis $\{E_1,E_2,E_5\}$ is a subalgebra of ${\g}$. Let $H$ be the Lie
subgroup of $G$ with a Lie algebra $b$. Having in mind $(24), (25)$ we
obtain that the section $\alpha =\{E_3, E_4\}$ is a normal to the submanifold
$H$, $\alpha $  is a totally real non-orthogonal to $\overline \xi $ section
of hybrid type and $\overline \xi \notin \alpha $, i.e. $H$ is of the same
type submanifolds considered in section 2.
We have the following decomposition of $\overline \xi , \overline {\varphi }X,
\overline {\varphi }E_3, \overline {\varphi }E_4$ with respect to $\{E_3, E_4\}$
and $T_pH$
\begin{equation}\label{26}
\begin {array}{llll}
\overline {\xi }=-\displaystyle{\frac{1}{2}}E_2+\displaystyle{\frac{\sqrt{3}}{2}}E_3; \\
\overline {\varphi }X=-\displaystyle{\frac{\sqrt{3}}{2}}\overline x^5E_2+
\displaystyle{\frac{\sqrt{3}}{2}}\overline x^2E_5-\displaystyle{\frac{1}{2}}
\overline x^5E_3+\overline x^1E_4; \\
\overline \varphi E_3=\displaystyle{\frac{1}{2}}E_5; \\
\overline \varphi E_4=-E_1;
\end{array} 
\end{equation}
where $X\in \chi(H)$ and
$X=\overline x^1E_1+\overline x^2E_2+\overline x^5E_5$.
We substitute 
\[
a=\displaystyle{\frac{\sqrt{3}}{2}}, \quad b=0, \quad k=\displaystyle{\frac{\sqrt{3}}{2}},
\quad 
\lambda =\displaystyle{\frac{4\sqrt{3}(2-\sqrt{3})}{3}}, \quad \mu =\lambda +1,
\]
\[
\xi _2=E_1, \quad \eta^2(X)=\overline x^1, \quad \varphi X=
-\displaystyle{\frac{\sqrt{3}}{2}}\overline x^5E_2+
\displaystyle{\frac{\sqrt{3}}{2}}\overline x^2E_5
\]
in $(10)$ and obtain an almost contact structure $(\phi ,\xi ,\eta )$ 
\begin{equation}\label{27}
\begin {array}{lll}
\xi =E_1; \\
\eta^(X)=\overline x^1; \\
\phi X=\displaystyle{\frac{2\sqrt{3}}{3}}\varphi X=-\overline x^5E_2+
\overline x^2E_5; \\
\end{array} 
\end{equation}
on the submanifold $H$.
\par
Using $(11), (16)$ and $(20)$ we get
\[
A_{E_3}X=-m\overline x^2E_1-m\overline x^1E_2; \quad
A_{E_4}X=-\displaystyle{\frac{1}{2}}\left(\displaystyle{\frac{3}{2}}a\overline x^2+
m\overline x^5\right)E_2+\displaystyle{\frac{1}{2}}m\overline x^2E_5;
\]
\[
\gamma (X)=\displaystyle{\frac{\sqrt{3}}{2}}\left(a\overline x^2-
m\overline x^5\right).
\]
Then the formulas of Gauss and Weingarten $(11)$ become 
\[
\begin {array}{lll}
\overline \nabla _XY=\nabla _XY-m\left(\overline x^2\overline y^1+\overline x^1\overline y^2\right)E_3+ 
\displaystyle{\frac{1}{2}}\left(\left(\displaystyle{\frac{3}{2}}a\overline x^2+
m\overline x^5\right)\overline y^2+m\overline x^2\overline y^5\right)E_4; \\ \\
\overline \nabla _XE_3=m\overline x^2E_1+m\overline x^1E_2+
\displaystyle{\frac{\sqrt{3}}{2}}\left(a\overline x^2-
m\overline x^5\right)E_4; \\ \\
\overline \nabla _XE_4=\displaystyle{\frac{1}{2}}\left(\displaystyle{\frac{3}{2}}a\overline x^2+
m\overline x^5\right)E_2-\displaystyle{\frac{1}{2}}m\overline x^2E_5+
\displaystyle{\frac{\sqrt{3}}{2}}\left(a\overline x^2-
m\overline x^5\right)E_3.
\end{array}
\]
Having in mind the last formulas, $(26)$ and $(27)$ we compute
the tensor $F$ of $H$
\[
F(X,Y,Z)=-\displaystyle{\frac{\sqrt{3}}{2}}a\overline x^2(\overline y^1\overline z^2+
\overline y^2\overline z^1), \quad X, Y, Z \in \chi (H)
\]
and verify that the submanifold 
$(H, \phi , \xi , \eta , g)$ 
of $G$, where $(\phi , \xi , \eta )$ is defined by (27)
is an almost contact manifold with Norden metric in the class 
${\F}_4\oplus {\F}_8$.

\end{document}